\input amstex
\documentstyle{amsppt}
\magnification=\magstep1
 \hsize 13cm \vsize 17.35cm \pageno=1
\loadbold \loadmsam
    \loadmsbm
    \UseAMSsymbols
\topmatter
%\NoRunningHeads
\title
On the analogs of Euler numbers and polynomials\\
 associated with  $p$-adic $q$-integral on $\Bbb Z_p$ at $q=-1$
\endtitle
\author
 Taekyun Kim
\endauthor
\affil
 {\it Jangjeon Research Institute for Mathematical Science and
 Physics,\\
 Ju-Kong Building 103-Dong 1001-Ho,\\
 Young-Chang Ri  544-4,  Hapchon-Gun Kyungnam,
 postal no. 678-802, S. Korea \\
 e-mail: tkim$\@$kongju.ac.kr, tkim64$\@$hanmail.net\\}
  \endaffil
\rightheadtext{ T. Kim  } \leftheadtext{$\lambda $-Euler numbers
and polynomials}

\abstract The purpose of this paper is to construct of
$\lambda$-Euler numbers and polynomials by using fermionic
expression of $p$-adic $q$-integral at $q=-1$. From these
$\lambda$-Euler polynomials, we derive the harmonic sums of higher
order. Finally, we investigate several interesting properties and
relationships involving the classical as well as the generalized
Euler numbers and polynomials. As an application, we will treat
$p$-adic invariant integrals on $\Bbb Z_p$ involving trigonometric
functions.
\endabstract
\thanks  2000 AMS Subject Classification: 11S80
\newline keywords and phrases : p-adic invariant integral,
trigonometric function, Euler numbers
\newline This paper was supported by Jangjeon Mathematical Society
(JMS-2006-R114006)
\endthanks
\endtopmatter

\document

{\bf\centerline {\S 1. Introduction}}

 \vskip 20pt

Let $p$ be a fixed odd prime, and let $\Bbb C_p$ be the $p$-adic
completion of the algebraic closure of $\Bbb Q_p$. For a fixed
positive integer $d$ with $(p,d)=1$, set
$$\split
& X=X_d = \lim_{\overleftarrow{N} } \Bbb Z/ dp^N \Bbb Z ,\cr & \
X_1 = \Bbb Z_p , \cr  & X^\ast = \underset {0<a<d p}\to {\cup} (a+
dp \Bbb Z_p ), \cr & a+d p^N \Bbb Z_p =\{ x\in X | x \equiv a
\pmod{dp^n}\},\endsplit$$ where $a\in \Bbb Z$ satisfies the
condition $0\leq a < d p^N$ [1].

The $p$-adic absolute value in $\Bbb C_p$ is normalized in such
way that $|p|_p =1/p$. Let $U_1 \subset \Bbb C_P$ denote the open
unit disc about $1$ and $U_d =\{ u\in \Bbb C_p | ~ |u^d -1|_p
<1\}$ the union of the open unit discs around $d$-th root of
unity. Let $U^m = U_d \times U_1^{m-1}$. When one talks of
$q$-extension, $q$ is variously considered as an indeterminate,  a
complex number $q\in\Bbb C$ or a $p$-adic number $q\in \Bbb C_p$.

If $q\in \Bbb C$, then we normally assume
$$|q-1|_p < p^{-1/p-1}\quad \text
{ so that }\quad  q^x =\exp (x \log q) \ \quad \text{ for }
|x|_p\leq 1.$$ Throughout this paper, we use the below notation :

$$[x]_q = [x:q]= \frac{1-q^x}{ 1-q} = 1+q +q^2 +\cdots +
q^{x-1}.$$

For $f \in C^{(1)} (\Bbb Z_p )(=$ the set of strictly
differentiable function on $\Bbb Z_p$), let us start with the
expressions
$$
\frac{1}{[p^N ]_q} \sum_{0\leq x < p^N } q^x f(x) = \sum_{0\leq x
<p^N } f(x) \mu_q (x +p^N \Bbb Z_p ),\quad  \text{cf. [11]},
$$
representing the $q$-analogue of the Riemann sums for  $f$. The
integral of $f$ on $\Bbb Z_p$ will be defined as the limit( $N \to
\infty) $ of these sums, when it exists. The $p$-adic $q$-integral
of a function $f\in C^{(1)}(\Bbb Z_p )$ is defined by author as
follows:
$$\eqalignno{
I_q (f) &=\int_{\Bbb Z_p} f(x) d \mu_q (x)&(1) \cr &=\lim_{N \to
\infty } \frac{1}{[p^N ]_q} \sum_{0\leq x < p^N} f(x) q^x , \text{
see [7,8,9,10].}&}
$$
For $f \in C^{(1)} (\Bbb Z_p )$, it is easy to see that
$$
\left| \int_{\Bbb Z_p } f(x) d \mu_q (x)  \right|_p \leq p ||f||_1
, \text{ cf. [2,3,4],}
$$
where $||f||_1 = \sup \Biggl\{ |f(0)|_p , \underset {x\neq y}\to
{\sup} \left| \dfrac{f(x) -f(y) }{ x-y}\right|_p \Biggr\}$.

If $f_n \to f$ in $C^{(1)} (\Bbb Z_p )$, namely $||f_n -f||_1 \to
0 $, then
$$
\int_{\Bbb Z_p} f_n (x) d\mu_q (x) \to \int_{\Bbb Z_p} f(x) d\mu_q
(x) , \text{ see [6,7,12]}.
$$

The $q$-analogue of the binomial coefficient is known as
$$
{{x} \brack {n}}_q =\dfrac{ [x]_q [x-1]_q \cdots [x-n+1]_q }{[n]_q
!},
$$
where $[n]_q ! = \prod_{i=1}^n [i]_q$, cf. [11].

 From this, we see
that
$$\eqalignno{ &
\int_{\Bbb Z_p} {{x} \brack {n}}_q d \mu_q (x) = \frac{(-1)^n
}{[n+1]_q} q^{n+1 - \binom{n+1}{2}}, \text{ (see [11]).}&(2)}
$$
For $q\in [0,1]$ certain $q$-deformed bosonic operators may be
introduced which generalize the undeformed bosonic ones (
corresponding $q=1$ ), see [13,14,15].

 The Eq. (2) is useful to
study $q$-deformed Stirling numbers $S_q (n,k)$ in the version
introduced by Milne [14]. From (1), we derive the usual ``bosonic"
$p$-adic $q$-integral as the limit $q\to 1$, i.e., $I_{q=1} (f)
=\int_{\Bbb Z_p} f(x) d\mu_{q=1} (x)$[4,5,6]. It is possible to
consider the case $q\in (-1,0)$ corresponding to $q$-deformed
fermionic certain and annihilation operators and the literature
given therein [16,22].

The expression (and the recursion relation ) for the $I_q (f)$
remains the same, so it is tempting to consider the limit $q\to
-1$. This limit will be called ``fermionic" since the
corresponding certain and annihilation operators are those of an
undeformed fermion. It is well known that the familiar Bernoulli
polynomials $B_n (x)$ are defined by means of the following
generating function:
$$\eqalignno{ &  F(x,t) = \frac{t}{e^t -1} e^{xt} =
\sum_{n=0}^\infty B_n (x) \frac{t^n}{n!}, \text{ cf.
[6,19,20],}&(3)}
$$
we note that, by substituting $x=0$ into (3), $B_n (0)=B_n$ are
the $n$-th Bernoulli numbers. By using ``bosonic" expression of
$p$-adic $q-$integral, one decade ago, we defined
$\lambda$-extension of these classical Bernoulli polynomials and
proved the properties of analogs to those satisfied by $B_n$ and
$B_n (x)$ [3,4,6,19]. That is,
$$
\eqalignno{  \lim_{q\to 1} \int_{\Bbb Z_p} f(x) \lambda^x d \mu_q
(x) &=\frac{\log \lambda +t}{\lambda e^t -1}\cr  &=
\sum_{n=0}^\infty B_n (\lambda ) \frac{t^n}{n!}.&(4)}
$$

Let $C_{p^n}$ be the cyclic group consisting of all $p^n -$th
roots of unity in $\Bbb C_p$ for any $n\geq 0$ and $T_p$ be the
direct limit of $C_{p^n}$ with respect to the natural morphisms,
hence $T_p$ is the union of all $C_{p^n}$ with discrete topology.

In the special $\lambda \in T_p$, we note that
$$\split
\lim_{q \to 1} \int_{\Bbb Z_p}  \lambda^x e^{tx} d \mu_q (x) &=
\frac{t}{\lambda e^t -1}\cr &= \frac{\log \lambda +t }{\lambda e^t
-1}\cr &=\sum_{n=0}^\infty B_n (\lambda) \frac{t^n }{n!},
\endsplit$$ where $B_n (\lambda)$ are $n $-th $\lambda $-Bernoulli
numbers [2,6]. Now, we consider the case $q\in (-1,0)$
corresponding to $q$-deformed fermionic certain and annihilation
operators and the literature given therein. The expression for the
$I_q (f)$ remains same, so it is tempting to consider the limit $q
\to -1$. That is,
$$\split
I_{q=-1} (f) &=\lim_{q\to -1} I_q (f) \cr &= \int_{\Bbb Z_p} f(x)
d \mu_{q=-1} (x), \text{ see [7,11].}\endsplit$$
 Thus, we note that

$$
\eqalignno{   \int_{\Bbb Z_p} e^{tx} \lambda^x d \mu_{q=-1} (x)
&=\frac{2}{\lambda e^t +1}\cr &= \sum_{n=0}^\infty E_n (\lambda )
\frac{t^n}{n!},&(5)}
$$
where $E_n (\lambda)$ are called $\lambda$-Euler numbers.

These numbers are classical and important in number theory. In
this paper we construct $\lambda$-Euler polynomials and numbers by
using ``fermionic" expression of $p$-adic $q$-integral. From these
$\lambda$-Euler polynomials, we derive the harmonic sums of higher
order.  Finally we investigate several interesting properties  and
relationships involving the classical as well as the generalized
Euler numbers and polynomials. As an application, we will treat
$p$-adic invariant integrals on $\Bbb Z_p$ involving trigonometric
functions.

\vskip 20pt

{\bf\centerline {\S 2. $\lambda$-Euler numbers and polynomials
associated }}

{\bf\centerline{with fermionic expression of $p$-adic
$q$-integral}}

\vskip 20pt

For $f\in C^{(1)} (\Bbb Z_p )$, we consider ``fermionic"
expression for the $p$-adic $q$-integral on  $\Bbb Z_p$ as
follows:
$$
\eqalignno{   I_{-1} (f) =I_{q=-1} (f) &= \lim_{q\to -1} I_q
(f)\cr &= \int_{\Bbb Z_p} f(x) d \mu_{-1} (x).&(6)}
$$

Let $f_1 (x)$ be translation with $f_1 (x) =f(x+1)$. Then we see
that
$$\split
I_{-1} (f_1 )& =- \lim_{n\to \infty} \sum_{x=0}^{p^N -1} f(x)
(-1)^x +2 f(0) \cr &= -I_{-1} (f) +2 f(0).\endsplit
$$
Therefore we obtain the following:

\proclaim{Lemma 1}  For $f\in C^{(1)}(\Bbb Z_p)$, we have
$$
\eqalignno{ &  I_{-1} (f_1 ) + I_{-1} (f) =2 f(0).&(7)}
$$
\endproclaim

From (7), we can easily derive the below theorem.

\proclaim{Theorem 2}  For $f\in C^{(1)}(\Bbb Z_p)$, $n\in \Bbb N$,
 we have
$$
\eqalignno{ &  I_{-1} (f_n ) + (-1)^n I_{-1} (f)
=2\sum_{x=0}^{n-1} (-1)^{n-1-x}f(x) ,&{}}
$$
where $f_n (x) =f(x+n)$.
\endproclaim

\vskip 10pt By Theorem 2, we can consider $\lambda$-Euler
polynomials. If we take $f(x)=\lambda^x e^{xt}\ (\lambda \in\Bbb
Z_p )$, then we have
$$\eqalignno{ &
\frac{2}{\lambda e^t +1} e^{x t} = \int_{\Bbb Z_p} e^{(x+y)t}
\lambda^y d \mu_{-1} (y). &(8)}$$

Now we define $\lambda$-Euler polynomials as follows:

$$\eqalignno{ &
F_\lambda (t,x ) =\frac{2}{\lambda e^t +1} e^{x t} =
\sum_{n=0}^\infty E_n (\lambda: x) \frac{t^n}{n!},&(9)}
$$
where $E_n (\lambda: x) $ are called $n$-th $\lambda$-Euler
polynomials. We also define, by substituting $x=0$ into (9), $E_n
(\lambda: 0)=E_n (\lambda)$ $n$-th $\lambda$-Euler numbers. From
(8) and (9), we derive the With's formula for $E_n (\lambda; x) $
as

$$\eqalignno{ & \int_{\Bbb Z_p} \lambda^y (x+y)^n d \mu_{-1} (y) =E_n (\lambda: x). &(9-1)}$$

Note that
$$\eqalignno{
E_n (\lambda: x)& =\sum_{l=0}^n \binom{n}{l} E_l (\lambda) x^{n-l}
\cr &=d^n \sum_{a=0}^{d-1}(-1)^a   \lambda^a E_n (\lambda^d :
\frac{a+x}{d}), &(10)}$$ where $d$ is positive integer. Because
$$\int_{\Bbb Z_p} (x+y)^n \lambda^y d \mu_{-1} (y)
=d^n \sum_{a=0}^{d-1} (-1)^a \lambda^a \int_{\Bbb Z_p} (
\frac{a+x}{d}+y )^n \lambda^{dy} d\mu_{-1} (y) .$$

\remark{Remark }  In [2,3,4,5,6,19], it was known that
$$\eqalignno{  \lim_{q\to 1} \int_{\Bbb Z_p} \lambda^y e^{(x+y)t}
d \mu_q (y)& = \frac{t+\log \lambda}{\lambda e^t -1} e^{xt} \cr &=
\sum_{m=0}^\infty B_m (\lambda ;x) \frac{t^m}{m!}, &(11)}$$ where
$B_m (\lambda ;x)$ are called $m$-th $\lambda$-Bernoulli
polynomials.

In view point of (11), we considered $\lambda$-Euler numbers and
polynomials.
\endremark

Let $\chi$ be the Dirichlet's character with conductor
$d(=\text{odd})\in\Bbb N$ and let us take $f(x) =\chi(x)
e^{tx}\lambda^x$.

From Theorem 2, we derive the below formula:
$$
\eqalignno{ &  \int_{\chi} \lambda^x e^{tx} \chi(x) d \mu_{-1} (x)
= \dfrac{2\sum_{a=0}^{d-1} e^{ta} (-1)^a \chi(a)
\lambda^a}{\lambda^d e^{dt} +1}. &(12)}
$$
By (12), we also consider the generalized $\lambda$-Euler numbers
attached to $\chi$ as follows:
$$\eqalignno{
F_{\lambda, \chi} (t) &=\dfrac{2\sum_{a=0}^{d-1} e^{ta} (-1)^a
\chi(a) \lambda^a}{\lambda^d e^{dt} +1}\cr & =\sum_{n=0}^{\infty}
E_{n,\chi} (\lambda) \frac{t^n}{ n! } . & (13)}
$$

From (12) and (13), we derive
$$\eqalignno{ &
\int_{\chi} \lambda^x  x^n  \chi(x) d \mu_{-1} (x) =
E_{n,\chi}(\lambda),\ \  n \geq 0.& (14)}
$$
It is easy to check that
$$\eqalignno{ &
\int_{\chi}  x^n  \chi(x) d \mu_{-1} (x) = d^n \sum_{a=0}^{d-1}
(-1)^a \lambda^a \chi(a)\int_{\Bbb Z_p}  \lambda^{xd}
\left(\dfrac{a}{d} +x \right)^n  d \mu_{-1} (x) .& (15)}
$$
From (9-1), (13) and (14), we derive
$$\eqalignno{ & E_{n,\chi}  (\lambda) = d^n
 \sum_{a=0}^{d-1}
(-1)^a \lambda^a \chi(a) E_n \left( \lambda^d :
\dfrac{a}{d}\right) .& (16)}
$$

\vskip 20pt

{\bf\centerline {\S 3. $\lambda$-zeta function associated with
$\lambda$-Euler numbers }}
\vskip 20pt

 In this section, we assume
that $\lambda(\neq -1)\in\Bbb C$ with $|\lambda |<1$. By (9), we
see that
$$\eqalignno{ & F_{\lambda} (t,x) =\frac{2}{\lambda e^t +1}e^{xt} =
2\sum_{n=0}^\infty (-1)^n  \lambda^n e^{(n+x)t} .& (17)}
$$
From (9) and (17), we note that
$$
E_k (\lambda : x) =\left. \frac{d^k}{ dt^k} F_\lambda (t,
x)\right\vert_{t=0} = 2 \sum_{n=0}^\infty (-1)^n \lambda^n (n+x)^k
.
$$
Therefore we obtain the following:

\proclaim{Theorem 3} For $k\in \Bbb N$, we have
$$
E_k (\lambda : x) = 2 \sum_{n=0}^\infty (-1)^n \lambda^n (n+x)^k.
$$
\endproclaim

\definition{Definition 4} For $s \in\Bbb C$, we define
$\lambda$-zeta function of Hurwitz's type as follows:
$$
\zeta_\lambda (s, x)= 2 \sum_{n=0}^\infty \frac{(-1)^n \lambda^n
}{(n+x)^s}.
$$
\enddefinition

Note that $\zeta_\lambda (s, x)$
 is analytic continuation in whole complex plane.  Let
$n$ be the even positive integer. Then we have
$$\eqalignno{
\dfrac{2(1- \lambda^n e^{nt})}{ \lambda e^t +1}&
=2\sum_{l=0}^{n-1} (-1)^l \lambda^l e^{lt} \cr &=
\sum_{m=0}^\infty (2 \sum_{l=0}^{n-1} (-1)^l \lambda^l l^m )
\frac{t^m}{m!}. & (18)}
$$

By (9) and (18), we see that
$$
\eqalignno{ & \sum_{m=0}^\infty (E_m (\lambda) -\lambda^n E_m
(\lambda: n)) \frac{t^m}{m!} = \sum_{m=0}^\infty
(2\sum_{l=0}^{n-1} (-1)^l \lambda^l l^m )\frac{t^m}{m!}. & (19)}
$$
By comparing the coefficients on the both sides in (19), we note
that
$$
\eqalignno{ & E_m (\lambda) -\lambda^n E_m (\lambda: n) =2
\sum_{l=0}^{n-1} (-1)^l \lambda^l l^m . & (20)}
$$
From (10) and (20), we derive
$$
\eqalignno{ &  2\sum_{l=0}^{n-1} (-1)^{l-1} \lambda^l l^m =
\lambda^n \sum_{l=0}^{m-1} \binom{m}{l} E_l (\lambda) n^{m-l} +
(\lambda^n -1 )E_m (\lambda) . & (21)}
$$
Therefore we obtain the following:

\proclaim{Theorem 5} Let $n$ be the positive even integer. Then we
have
$$
\eqalignno{ &  2\sum_{l=0}^{n-1} (-1)^{l-1} \lambda^l l^m =
\lambda^n \sum_{l=0}^{m-1} \binom{m}{l} E_l (\lambda) n^{m-l} +
(\lambda^n -1 )E_m (\lambda) . & }
$$
\endproclaim
\remark{Remark}
$$
\eqalignno{ & E_0 (\lambda) = \frac{2}{\lambda +1},\   E_1
(\lambda) =- \frac{2 \lambda}{ (\lambda+1)^2},\  E_2 (\lambda)
=\frac{4\lambda^2 -2 \lambda +2 }{ (\lambda+1)^3} ,\cdots .&(22)}
$$
\endremark

\proclaim{Corollary 6} Let $k$ be  positive  integer. Then we see
that
$$ \zeta_\lambda (-k, x ) = E_k (\lambda : x).
$$
\endproclaim

Let $\chi$ be the primitive Dirichlet's character with conductor
$d$(=odd)$\in \Bbb N$. From (13), we derive
$$
\eqalignno{  F_{\lambda, \chi} (t) &= \frac{2 \sum_{a=0}^{d-1}
e^{ta} (-1)^a \lambda^a \chi (a)}{\lambda^d e^{dt} +1}\cr  =& 2
\sum_{n=0}^\infty (-1)^n \chi (n) \lambda^n e^{nt}.&(23)}
$$

By (13) and (23), we easily see that
$$\split E_{k,\chi}(\lambda) &=\left. \frac{d^k}{ dt^k} F_{\lambda , \chi}(t)
\right\vert_{t=0} \cr &= 2 \sum_{n=1}^\infty \chi (n) \lambda^n
n^k .
\endsplit
$$
Therefore we obtain the following:

\proclaim{Theorem 7} Let $\chi$ be  the primitive Dirichlet's
character with conductor $d$(=odd)$\in \Bbb N$. Then, we have
$$ E_{k,\chi} (\lambda)  =2\sum_{n=1}^\infty (-1)^n \chi (n) n^k \lambda^n.
$$
\endproclaim

\definition{Definition 8} For $s\in \Bbb C$, we define
$\lambda-l-$function as follows:
$$
\eqalignno{ & l_{\lambda}(s, \chi) = 2  \sum_{n=1}^{\infty}
\dfrac{(-1)^n \chi(n) \lambda^n}{n^s} .&(24)}
$$
\enddefinition
Note that $l_{\lambda}(s, \chi)$ is also analytic continuation in
whole complex plane.

\proclaim{Corollary 9} Let $k\in \Bbb N$. Then we have
$$
l_\lambda (-k, \chi ) = E_{k,\chi} (\lambda).
$$
\endproclaim

Let $s$ be a complex variable, and let $a$ and $F($=odd) be
integer with $0< a< F$. We consider the below harmonic sum(
partial $\lambda$-zeta function):
$$\split
H_\lambda (s,a|F)&=\sum_{ { m\equiv a (F)} \atop {m>0} } \dfrac{
(-1)^m \lambda^m }{ m^s} \cr &= \sum_{m=0}^\infty
\dfrac{\lambda^{nF +a} (-1)^{nF +a} }{ (a+nF)^s} \cr &= \lambda^a
(-1)^a \sum_{n=0}^\infty \dfrac{ (\lambda^F )^n (-1)^n }{ (a+
nF)^s } \cr &= \dfrac{\lambda ^a (-1)^a  }{2} F^{-s}
\zeta_{\lambda^F} (s, \dfrac{a}{F} ).\endsplit
$$

Note that

$$
\eqalignno{ & H_{\lambda}(-k, a|F) = \dfrac{(-1)^a \lambda^a
F^k}{2}E_k (\lambda^F : \dfrac{a}{F} ).&(25)}
$$
From (24) and (25), we derive
$$
l_\lambda (s,\chi ) =2 \sum_{a=1}^F \chi (a) H_\lambda (s, a|F).
$$
The harmonic sum $H_\lambda (s, a|F )$ will be called partial
$\lambda-$zeta function which interpolates $\lambda$-Euler
polynomials at negative integers [21].  The values of $l_\lambda
(s,\chi)$ at negative integers are algebraic, hence may be
regarded as lying in an extension of $\Bbb Q_p$. We therefore look
for a $p$-adic function which agrees with $l_\lambda (s,\chi)$ at
the negative integers in  Section 4. \vskip 20pt

{\bf\centerline {\S 4. $p$-adic harmonic sums of higher order
associated with $p$-adic $\lambda-l-$function}} \vskip 20pt

Let $w(x)$ be the Teichm\"uller character  and let $<x>
=\dfrac{x}{w(x)}$. When $F$(=odd) is a multiple of $p$ and
$(a,p)=1$, we define $p-$adic partial $\lambda-$zeta function as
follows :

$$
H_{\lambda, p} (s, a|F) =\dfrac{(-1)^a \lambda^a}{2}<a>^{-s}
\sum_{j=0}^\infty \binom{-s}{j} \left( \dfrac{F}{a}\right)^j E_j
(\lambda^F ),
$$
for $s\in \Bbb Z_p$. Thus, we note that
$$
\split H_{\lambda, p} (-n , a|F) & = \frac{(-1)^a \lambda^a}{2}
<a>^n \sum_{j=0}^n \binom{n}{j} \left( \dfrac{F}{a}\right)^j E_j
(\lambda^F )\cr &= F^n  \dfrac{(-1)^a \lambda^a }{2} w^{-n} (a)
\sum_{j=0}^n \binom{n}{j}\left( \dfrac{a}{F}\right)^{n-j} E_j
(\lambda^F)\cr
 &= \dfrac{(-1)^a}{2} F^n w^{-n} (a) E_n (\lambda^F :
 \frac{a}{F})\cr
 &= w^{-n}(a) H_\lambda (-n, a|F),
\endsplit
$$
where $n$ is positive integer.

Now we consider $p-$adic interpolating function for
$\lambda-$Euler numbers as follows:
$$
l_{\lambda, p} (s, \chi) = 2 \sum_{{a=1} \atop {(a,p)=1}}^F
 \chi(a) H_{\lambda, p}(s, a|F),$$
 for $s\in \Bbb Z_p$. Let $n$ be natural numbers. Then we have
 $$\split l_{\lambda, p} (-n, \chi)& =2 \sum_{{a=1} \atop {(a,p)=1}}^F
 \chi(a) H_{\lambda, p}(-n, a|F)\cr
 & =E_{n,\chi w^{-n}}(\lambda) -p^n
 \chi w^{-n}(p) E_{n,\chi w^{-n}} (\lambda^p ).\endsplit
 $$
In fact, we see that
$$
l_{\lambda, p} (s, \chi) =\sum_{a=1}^F (-1)^a \lambda^a \chi (a)
\sum_{j=0}^\infty \binom{-s}{j} \left( \dfrac{F}{a} \right)^j E_j
(\lambda),
$$
for $s\in \Bbb Z_p$. This is a $p-$adic analytic function and has
the following properties for $\chi =w^t$:
$$
\eqalignno{ & l_{\lambda, p}(-n, w^t ) =E_n (x) -p^n E_n
(\lambda^p ),&(26)}
$$
where $n\equiv t \pmod{p-1}$, $l_{\lambda, p}(s, w^t )\in\Bbb Z_p$
for all $s\in\Bbb Z_p$  when $t\equiv 0\pmod{p-1}$.

If $t\equiv 0\pmod{p-1}$, then $l_{\lambda, p}(s, w^t )\equiv
l_{\lambda, p}(s_2 , w^t )\pmod{p}$ for all $s,s_2 \in\Bbb Z_p$.

It is easy to see that
$$
\dfrac{1}{ r+k -1} \binom{-r}{k} \binom{1-r-k}{j} = \dfrac{-1}{
j+k} \binom{-r}{k+j -1} \binom{k+j}{j},
$$
for all positive integers $r,j,k$ with $j,k\geq 0$, $j+k>0$, and
$r\neq 1-k$.

Thus, we note that
$$
\dfrac{r}{ r+k } \binom{-r-1}{k} \binom{-r-k}{j} = \binom{-r}{k+j
} \binom{k+j}{j}.
$$
For $F$(=odd)$\in\Bbb N$, let $n$ be positive even integer then we
have
$$
\eqalignno{  \sum_{l=0}^{n-1} \dfrac{\lambda^{Fl+a} (-1)^{Fl +a}
}{ (Fl+a)^r } & = \sum_{l=0}^{n-1}  (-1)^{Fl +a} a^{-r}
\lambda^{Fl +a}  \sum_{s=0}^{\infty } \binom{-r}{s} \left(
\dfrac{Fl }{a}\right)^s \cr &= \sum_{m=0}^{\infty } \binom{-r}{m}
a^{-r} (-1)^a \lambda^a  \left( \dfrac{F}{a}  \right)^m
\sum_{l=0}^{n-1} (-1)^l l^m \lambda^{Fl} &{(27)}\cr
 &=- \sum_{m=0}^{\infty }  \binom{-r}{m} a^{-r} (-1)^a \lambda^a \left(
 \dfrac{F}{a} \right)^m \dfrac{1}{2} \{ \lambda^{Fn} \sum_{l=0}^{m-1} \binom{m}{l}
 E_l (\lambda^F ) n^{m-l}\cr &\quad +
 (\lambda^{Fn }-1)E_m (\lambda^F )\}.
 }
$$

Hence, we note that
$$\eqalignno{ &
-\sum_{l=0}^{n-1} \dfrac{\lambda^{Fl+a} (-1)^{Fl +a} }{ (Fl+a)^r }
- \dfrac{(\lambda^{Fn} -1 )}{2} \sum_{m=0}^\infty \binom{-r}{m}
a^{-r} (-1)^a
 \lambda^a  \left( \dfrac{F}{a}  \right)^m E_m (\lambda^F)\cr
 &=\dfrac{1}{2} \sum_{m=0}^\infty \binom{-r}{m}a^{-r} (-1)^a
 \lambda^a \left( \dfrac{F}{a}  \right)^m \lambda^{Fm}
 \sum_{l=0}^{m-1} \binom{m}{l} n^{m-l} E_l (\lambda^F )\cr
 &=\lambda^{Fn} \dfrac{(-1)^a a^{-r} \lambda^a}{2}
 \sum_{s=0}^\infty\binom{-r}{s} \left( \dfrac{F}{a}  \right)^s
 \left( \sum_{l=0}^{s-1} \binom{s}{l} n^{s-l} E_l (\lambda^F )
 \right)\cr
 &=\lambda^{Fn} \sum_{s=0}^\infty \binom{-r}{s} w^{-r} (a) \left( \dfrac{F}{a}
 \right)^s \dfrac{(-1)^a \lambda^a}{2} <a>^{-r} \sum_{l=0}^{s-1}
 \binom{s}{l} n^{s-l} E_l (\lambda^F)\cr
 &=\lambda^{Fn} \sum_{k=0}^\infty\sum_{l=0}^\infty \binom{-r}{k+l} w^{-r} (a)
 \left( \dfrac{a}{F} \right)^{-k-l} n^k \dfrac{(-1)^a
 \lambda^a}{2}<a>^{-r} E_l (\lambda^F ) \binom{k+l}{l}\cr
 &=\lambda^{Fn} \sum_{k=0}^\infty\sum_{l=0}^\infty  \dfrac{r}{r+k}
 \binom{-r-1}{k} \binom{-r-k}{l} a^{-r} \left( \dfrac{a}{F} \right)^{-k-l}
 n^k
 \frac{(-1)^a}{2}\lambda^a E_l (\lambda^F )\cr
& = \lambda^{Fn} \sum_{k=0}^\infty \dfrac{r}{r+k} \binom{-r-1}{k}
w^{-k-r} (a) (nF)^k \frac{(-1)^a}{2}\lambda^a <a>^{-k-r}
\sum_{l=0}^\infty \binom{-r-k}{l} \left( \dfrac{F}{a}\right)^l E_l
(\lambda^F )\cr
 &= \lambda^{Fn}  \sum_{k=0}^\infty \dfrac{r}{r+k}  \binom{-r-1}{k} w^{-k -r}(a) (nF)^k
 H_{\lambda, p}(r+k, a|F).   &(28)}
$$
For $F=p$, $r\in\Bbb N$, $n$(=even)$\in\Bbb N$, we have
$$\eqalignno{ &
2\sum_{{j=1} \atop {(j,p)=1}}^{np} \dfrac{(-1)^j \lambda^j }{j^r}
=2 \sum_{a=1}^{p-1} \sum_{l=0}^{n-1} \dfrac{(-1)^{a+pl}
\lambda^{a+pl}}{(a+pl)^r}.&(29)}
$$
Now, we set
$$\eqalignno{ &
B^{(r)} (a,F) = \dfrac{1}{2} \sum_{m=0}^\infty \binom{-r}{m}a^{-r}
(-1)^a \lambda^a \left(  \dfrac{F}{a} \right)^m E_m (\lambda^F
).&(30)}
$$
From (27), (28), (29) and (30), we derive
$$\split &
 2\sum_{{j=1} \atop {(j,p)=1}}^{np} \dfrac{(-1)^j \lambda^j
}{j^r}\cr &=2 \sum_{a=1}^{p-1} \sum_{l=0}^{n-1} \dfrac{(-1)^{a+pl}
\lambda^{a+pl}}{(a+pl)^r}\cr
 &= -2\sum_{a=1}^{p-1} \left\lbrace \lambda^{pn}
\sum_{k=0}^\infty \dfrac{r}{r+k}  \binom{-r-1}{k} w^{-k-r} (a)
(np)^k H_{\lambda, p} (r+k, a|F) \right. \cr &\quad \left.
+(\lambda^{pn} -1)B^{(r)}(a,p) \right\rbrace\cr
&=-\sum_{k=0}^\infty \dfrac{r}{r+k} \binom{-r-1}{k} \lambda^{pn}
(pn)^k  l_{\lambda,p}(r+k, w^{-k-r}) -2 (\lambda^{pn}
-1)\sum_{a=1}^{p-1} B^{(r)} (a,p).
\endsplit
$$
Therefore we obtain the following :

\proclaim{Theorem 10} Let $p$ be an odd prime and let $n$ be  an
even positive integer. If $r\in\Bbb N$, then we have
$$\split
2\sum_{{j=1} \atop {(j,p)=1}}^{np} \dfrac{(-1)^j \lambda^j }{j^r}
&= -\sum_{k=0}^\infty \dfrac{r}{r+k} \binom{-r-1}{k} \lambda^{pn}
(pn)^k l_{\lambda,p} (r+k, w^{-k-r})\cr &\quad -2 (\lambda^{pn}
-1) \sum_{a=1}^{p-1} B^{(r)} (a,p) .\endsplit$$
\endproclaim

\proclaim{Corollary  11} Let $\lambda=1$. Then we see [21] that
$$\split
2\sum_{{j=1} \atop {(j,p)=1}}^{np} \dfrac{(-1)^j }{j^r} &=
-\sum_{k=0}^\infty \dfrac{r}{r+k} \binom{-r-1}{k} (pn)^k l_{1,p}
(r+k, w^{-k-r}) .\endsplit$$
\endproclaim

\remark{Remark} Luo-Srivastava have studied $\lambda-$Bernoulli
numbers and polynomials [17, 18]. In [17, 18], $\lambda-$Bernoulli
numbers and polynomials are called by Apostol-Bernoulli numbers
and polynomials. Indeed, these numbers are Euler numbers[2,3,6].
In [6], we have studied these numbers and polynomials by using the
``bosonic" expression of $p-$adic $q-$integral on $\Bbb Z_p$. That
is, $\dfrac{ \log \lambda +t }{ \lambda e^t -1} =\sum_{n=0}^\infty
B_n (\lambda) \dfrac{t^n}{n!}$.

However, this generating function is not simple. So, we modified
these generating function of $B_n (\lambda)$ in the locally
constant space.

If $\lambda \in T_p= \underset {n\geq 0} \to {\cup} \{ \zeta |
\zeta^{p^n}=1\}$, then we have
$$
\dfrac{ \log \lambda +t}{ \lambda e^t -1} =\dfrac{t}{ \lambda e^t
-1} =\sum_{n=0}^\infty B_n (\lambda) \dfrac{t^n}{n!}, \ \
\text{cf. } [2,3,4,6,19].
$$
Recently, Simsek, Jang-Pak-Rim also studied for the properties of
$\lambda-$Bernoulli numbers and polynomials(see [2,3,19]). In
[16], M. Schork also studied the representation of the
$q-$fermionic commutation relations and limit $q=1$.
\endremark

\vskip 20pt

{\bf\centerline {\S 5. Applications of $p$-adic $q$-integrals at
$q=-1$ (or $q=1$) }} \vskip 20pt

In this section we consider $p$-adic $q$-integral at $q=-1$ (or
$q=1$ ) involving trigonometric functions. By Eq.(7), we see that
$ I_{-1} (f_1 ) + I_{-1} (f) =2 f(0)$, for $f\in C^{(1)}(\Bbb
Z_p)$. If we take $f(x)= \cos ax $ ( or $f(x)=\sin ax $), then we
have
$$\aligned
0& =\int_{\Bbb Z_p}\sin a(x+1) d\mu_{-1}(x)+\int_{\Bbb Z_p}\sin
ax d\mu_{-1}(x)\\
 &= (\cos a+1)\int_{\Bbb Z_p} \cos ax
d\mu_{-1}(x)+\sin a\int_{\Bbb Z_p}\cos ax d\mu_{-1}(x),
\endaligned\tag 31$$
and
$$\aligned
2& =\int_{\Bbb Z_p}\cos a(x+1) d\mu_{-1}(x)+\int_{\Bbb Z_p}\cos
ax d\mu_{-1}(x)\\
 &= (\cos a+1)\int_{\Bbb Z_p} \cos ax
d\mu_{-1}(x)-\sin a\int_{\Bbb Z_p}\sin ax d\mu_{-1}(x).
\endaligned \tag32$$
Thus, we note tahat

$$\aligned
&\int_{\Bbb Z_p}\sin ax d\mu_{-1}(x)=-\frac{\sin a}{\cos a+1}
\int_{\Bbb Z_p}\cos ax d\mu_{-1}(x),\\
 & 2(\cos a +1)=(\cos a+1)^2\int_{\Bbb Z_p} \cos ax
d\mu_{-1}(x)+\sin ^2a\int_{\Bbb Z_p}\cos ax d\mu_{-1}(x).
\endaligned \tag33$$
From (31), (32) and (33), we derive
$$
\int_{\Bbb Z_p} \cos ax d\mu_{-1}(x)=1, \int_{\Bbb Z_p} \sin ax
d\mu_{-1}(x)=-\frac{\sin a}{\cos a+1}. \tag 34$$ Therefore we
obtain the below:

\proclaim{ Proposition 12} Let $ a\in \Bbb Z_p$. Then we have
$$
\int_{\Bbb Z_p} \cos ax d\mu_{-1}(x)=1, \int_{\Bbb Z_p} \sin ax
d\mu_{-1}(x)=-\frac{\sin a}{\cos a+1}=-\tan \frac{a}{2}. $$
\endproclaim

By using Taylor expansion for $\sin ax$ at $x=0$,  we easily see
that
$$ \sin
ax=\sum_{n=0}^{\infty}\frac{(-1)^n}{(2n+1)!}a^{2n+1}x^{2n+1}.
\tag35$$ From (9-1), (34) and (35), we derive
$$-\tan \frac {a}{2}=\int_{\Bbb Z_p} \sin ax
d\mu_{-1}(x)=\sum_{n=0}^{\infty}\frac{(-1)^na^{2n+1}}{(2n+1)!}\int_{\Bbb
Z_p}x^{2n+1}d\mu_{-1}(x).\tag 36$$

\proclaim{Theorem 13 } Let $\lim_{q\rightarrow -1}\int_{\Bbb
Z_p}f(x)d\mu_q(x)=\int_{\Bbb Z_p}f(x)d\mu_{-1}(x).$ Then we have
$$\tan \frac {a}{2}=-\int_{\Bbb Z_p} \sin ax
d\mu_{-1}(x)=\sum_{n=0}^{\infty}\frac{(-1)^{n+1}a^{2n+1}}{(2n+1)!}E_{2n+1},$$
where $E_{n}$ are the $n$-th Euler numbers.
\endproclaim
 Let $f_1 (x)=f(x+1)$ be translation in $C^{1}(\Bbb Z_p)$ and let
 $f^{\prime} (0)=\frac{d}{dx}f(x)|_{x=0}$.  Then we note that
  $$I_1 (f)=\lim_{q\rightarrow 1}\int_{\Bbb Z_p} f(x)
d\mu_{q}(x )\Longrightarrow I_{1}(f_1)-I_{1}(f)=f^{\prime}(0).
\tag 37$$ By (37), we  easily see that
$$\aligned
0&=\int_{\Bbb Z_p}\cos a(x+1)d\mu_1(x)-\int_{\Bbb Z_p} \cos ax
d\mu_1(x)\\
&=(\cos a -1)\int_{\Bbb Z_p} \cos ax d\mu_1(x)-\sin a\int_{\Bbb
Z_p}\sin ax d\mu_1(x).
\endaligned$$
Thus, we have
$$\int_{\Bbb Z_p} \sin ax d\mu_1(x)=-\frac{a}{2}. \tag 38$$
It is easy to see that
$$\aligned
 a&=\int_{\Bbb Z_p}\sin a(x+1)d\mu_1(x)-\int_{\Bbb Z_p}\sin ax
 d\mu_1(x)\\
 &=(\cos a -1)\int_{\Bbb Z_p}\sin ax d\mu_1(x)+\sin a\int_{\Bbb
 Z_p}\cos ax d\mu_1(x).
\endaligned\tag39$$
From (38) and (39), we derive
$$\int_{\Bbb Z_p} \cos ax d\mu_1(x)=\frac{a\sin a}{2-2\cos
a}=\frac{a}{2}cot \frac{a}{2}. \tag 40$$ By using Taylor
expansion, we easily see that
$$\cos ax =\sum_{n=0}^{\infty}\frac{(-1)^na^{2n}}{(2n)!}x^2.$$
Thus, we have
$$\frac{a}{2}\cot a=\int_{\Bbb Z_p} \cos ax
d\mu_1(x)=\sum_{n=0}^{\infty}\frac{(-1)^nB_{2n}}{(2n)!}a^2n ,$$
where $B_n$ are the $n$-th Bernoulli numbers.

\enddocument